\documentclass[12pt,reqno]{article}

\usepackage[english]{babel}
\usepackage{amsmath}
\usepackage{amssymb}
\usepackage{amsthm}

\usepackage{fullpage}
\usepackage[usenames]{color}
\usepackage{graphicx}
\usepackage[colorlinks=true,
linkcolor=webgreen,
filecolor=webbrown,
citecolor=webgreen]{hyperref}

\definecolor{webgreen}{rgb}{0,.5,0}
\definecolor{webbrown}{rgb}{.6,0,0}
\definecolor{red}{rgb}{1,0,0}

\usepackage{color}

\newcommand{\seqnum}[1]{\href{http://www.research.att.com/cgi-bin/access.cgi/as/~njas/sequences/eisA.cgi?Anum=#1}{\underline{#1}}}

\newcommand{\mytext}[1]{ \: \textrm{#1} \: }
\newcommand{\mysetdescr}[2]{\left\{ #1 \: \left| \: #2 \right. \right\} }
\newcommand{\streleq}[1] {\stackrel{\eqref{#1}}{=}}

\newcommand{\darr}{{\downarrow \,}}
\newcommand{\uarr}{{\uparrow \,}}
\newcommand{\udarr}{{\updownarrow \,}}

\newcommand{\setx}[1]{ \{ #1 \} }
\newcommand{\myN}{\mathbb{N}}

\newcommand\urbild[1]{#1^{-1}}

\def\D{{\cal D}}
\def\H{{\cal H}}
\def\R{{\cal R}}

\newcommand{\mf}[1]{\mathfrak{ #1 }}
\newcommand{\fT}{\mf{T}}

\def\BP{\begin{proof}}
\def\EP{\end{proof}}

\newcommand{\Ch}[1]{{\bf #1}}
\newcommand{\ACh}[1]{\overline{\Ch{#1}}}

\newcommand{\Brd}{B_{-}}
\newcommand{\Bru}{B^{-}}
\newcommand{\Brr}{B_{-}^{-}}

\newcommand{\brd}{b_{-}}
\newcommand{\bru}{b^{-}}
\newcommand{\brr}{b_{-}^{-}}

\newcommand{\binomk}[2]{{#1 \choose #2 }}

\begin{document}

\theoremstyle{plain}

\newtheorem{theorem}{Theorem}
\newtheorem{corollary}{Corollary}
\newtheorem{lemma}{Lemma}

\theoremstyle{definition}

\newtheorem{definition}{Definition}

\title{\bf A Flexible Approach for the Enumeration of Down-Sets and its Application on Dedekind Numbers}
\author{Frank a Campo}
\date{\small 41747 Viersen, Germany\\
{\sf acampo.frank@gmail.com}}

\maketitle

\begin{abstract}
We introduce a flexible approach for the enumeration of the down-sets of a finite poset and test it with the calculation of the Dedekind numbers $b(5) = 7581$ and $b(6) = 7828354$. For the calculation of $b(5)$, we develop two methods of which the first one (without pre-calculations) requires simple evaluation of 80 posets and the second one (with pre-calculations) of 34 posets. The calculation of $b(6)$ (with pre-calculations) is done by evaluating 245 posets.

\end{abstract}

\section{Introduction}

The enumeration of the down-sets of a finite poset is an important task in order combinatorics, but unfortunately, it quickly becomes complicated. For simple posets and simple ways to generate them (e.g., chains, antichains, direct and ordinal sums), folklore formulas exist for the number of down-sets, but even for a standard structure like the product of posets, no handy general formula is available. For an overview about the enumeration of down-sets, the reader is referred to literature \cite{aCampo_2018_Framework,aCampo_Erne_2019}.

For a poset $P$, let $\D(P)$ denote its down-set lattice and $d(P) := \# \D(P)$ the number of its down-sets. A common approach for the calculation of a parameter of a complicated structure is to calculate it by evaluating substructures, in our case, to find sub-posets $S_1, \ldots ,S_n$ of $P$ with $d(P) = d(S_1) + \cdots + d(S_n)$. A textbook example \cite[Ex.\ 1.14]{Davey_Priestley_2012} is the formula
$$
d(P) = d(P \setminus{x}) + d(P \setminus ( \darr x \cup \uarr x ) )
$$
with $x$ being an arbitrary point of $P$. An advanced result has been proven by the author and Ern\'{e} \cite[Theorem 3.1]{aCampo_Erne_2019} in which the summation runs over all subsets of an arbitrary {\em  antichain} in $P$. In Theorem \ref{theo_dP_formel} in Section \ref{sec_countingDownsets}, we present a generalization of this result: For an arbitrary  {\em subset} $M$ of the carrier of $P$, we calculate $d(P)$ by running a summation over all down-sets of the poset induced by $P$ on $M$.

The power of Theorem \ref{theo_dP_formel} lies in its flexibility. For the enumeration of the down-sets of $P$, we can select the set $M$ in such a way that the resulting summation process is simple in the following sense:
\begin{enumerate} \label{Kriterien einfachheit}
\item With $Q$ being the poset induced by $P$ on $M$, the number of down-sets of $Q$ is small and they can easily be generated.
\item For each down-set $D \in \D(Q)$, it is easy to calculate the respective summand $d(S_D)$.
\item All steps are easy to programme.
\end{enumerate}

\begin{table}
\begin{center}
\begin{tabular}{| c | r | }
\hline
$ n $ & $n$\textsuperscript{th} Dedekind number \\
\hline \hline
0 &                       2 \\
1 &                       3 \\
2 &                       6 \\
3 &                      20 \\
4 &                     168 \\
5 &                    7581 \\
6 &                 7828354 \\
7 &           2414682040998 \\
8 & 56130437228687557907788 \\
\hline
\end{tabular}
\caption{\label{table_Dedekind_Zahlen} Known Dedekind numbers.}
\end{center}
\end{table}

In order to test our approach, we apply it on the calculation of Dedekind numbers. (The sequence number of Dedekind numbers in OEIS \cite{Sloane_DB} is \seqnum{A000372}.) These numbers have attracted the interest of many scientists because they are so difficult to calculate. It is possible to reduce the effort considerably by performing pre-calculations, but also these are challenging. In fact, only the first nine Dedekind numbers are known (Table \ref{table_Dedekind_Zahlen}). For references, see Section \ref{sec_BooleanLatt}; more details and a historical overview are found in literature \cite{aCampo_2018}.

The $n$\textsuperscript{th} Dedekind number is the number of down-sets (up-sets, antichains) of the Boolean lattice $B(n)$ with $n$ atoms. By applying Theorem \ref{theo_dP_formel} on $B(n)$, we get in Theorem \ref{theo_Dk} in Section \ref{sec_BooleanLatt} a formula for the $n$\textsuperscript{th} Dedekind number $b(n) := d(B(n))$ in which all but one parameter can recursively be determined. The exception is the parameter $\brr(n)$, the number of down-sets of the poset $\Brr(n)$ we get by removing both extrema, all atoms, and all co-atoms from $B(n)$.

In Section \ref{sec_Application}, we develop different methods to calculate $\brr(5) = 6212$ and $\brr(6) = 7741776$ by applying Theorem \ref{theo_dP_formel} in different ways on $\Brr(5)$ and $\Brr(6)$. The reason for us to restrict us to the moderate cases $n = 5$ and $n = 6$ is that we have to do all calculations with pencil, paper, and ordinary table calculation.

In the calculation of $\brr(5)$, our best methods require a summation over 80 simple evaluations of down-sets without any pre-calculation and over 34 evaluations if we invest in the pre-calculation of isomorphism classes. In the calculation of $\brr(6)$, we pre-calculate a table with 1024 entries, and using isomorphism classes, we have to evaluate 245 down-sets. 

In the application, our intention was to test the flexibility of our approach to generate efficient methods for the enumeration of down-sets, and we do not claim that we have found the most efficient methods for the calculation of $b(5)$ and $b(6)$. Our numbers of down-set evaluations are clearly better than those required by the standard algorithm for the calculation of Dedekind numbers in its simple form (cf.\ Section \ref{sec_BooleanLatt}), but we cannot use this algorithm as benchmark because our 34 and 245 evaluations have been achieved with taking isomorphism into account. (Moreover, the algorithm tested by Fidytek et al.\ \cite{Fidytek_etal_2001} is faster than the standard algorithm.) The performance parameters of algorithms exploiting isomorphism reported in literature \cite{Fidytek_etal_2001,Markowsky_1989} are based on run-time of program execution in the calculation of $b(7)$. But we did not calculate $b(7)$, and run-time was not relevant in our work. A comparison of our methods with the algorithms found in literature is thus not possible at the present state. 

\section{Notation} \label{sec_notation}

We are working with {\em finite partially ordered sets (posets)}, that is ordered pairs $P = (X,\leq_P)$ consisting of a finite set $X$ (the {\em carrier} of $P$) and a {\em partial order relation} $\leq_P$ on $X$, i.e., a reflexive, antisymmetric, and transitive subset of $X  \times X$. We define
\begin{align*}
{<_P} & := \mysetdescr{ (x,y) \in {\leq_P} }{ x \not= y },
\end{align*}
and as usual, we write $x \leq_P y$ and $u <_P v$ for $(x,y) \in {\leq_P}$ and $(u,v) \in {<_P}$. For $x <_P y$, we say that the point $y$ {\em covers} the point $x$, iff there exists no $z \in X$ with $x <_P z <_P y$.

For $Y \subseteq X$, the {\em induced sub-poset} $P \vert_Y$ of $P$ is $\left( Y, {\leq_P} \cap (Y \times Y) \right)$. To simplify notation, we identify a subset $Y \subseteq X$ with the poset $P \vert_Y$ induced by it. Furthermore, we write $P - Y$ instead of $P \vert_{X \setminus Y}$.

We call a poset $P = (X, \leq_P)$ an {\em antichain} iff its partial order relation is the {\em diagonal} $\mysetdescr{ (x,x) }{ x \in X }$, and we call it a {\em chain} iff $x \leq_P y$ or $y \leq_P x$ holds for all $x, y \in X$. We write $\ACh{a}$ for an antichain with $a$  points and $\Ch{c}$ for a chain with $c$ points. For the sake of simplicity, we assume that $\Ch{c}$ has the carrier $0, \ldots ,c-1 $ equipped with the natural order.

A subset $Y \subseteq X$ is called a {\em down-set (up-set)} of $P$, iff $x \leq_P y$ implies $x \in Y$ for all $y \in Y$ and $x \in X$ (iff $y \leq_P x$ implies $x \in Y$ for all $y \in Y$ and $x \in X$). For $y \in X$, we define the {\em down-set} and {\em up-set induced by $y$} as
\begin{align*}
\darr_P \; y & := \mysetdescr{ x \in X }{ x \leq_P y }, \\
\uarr_P \; y & := \mysetdescr{ x \in X }{ y \leq_P x },
\end{align*}
and for $V \subseteq X$, we define $\darr_P V := \cup_{ v \in V } \darr_P v$ and $\uarr_P V := \cup_{ v \in V } \uarr_P v$. The symbol $\D(P)$ denotes the set of down-sets of $P$, and $d(P) := \# \D(P)$ is the cardinality of $\D(P)$. Together with set-inclusion, $\D(P)$ is a lattice.

For posets $P = (X, \leq_P)$ and $Q = (Y, \leq_Q)$, we define their {\em product} $P \times Q = ( X \times Y, \leq_{P \times Q} )$ by
$$
(x_1, y_1) \leq_{P \times Q} (x_2, y_2) \quad \Leftrightarrow \quad x_1 \leq_P x_2 \; \mytext{and } y_1 \leq_Q y_2.
$$
The symbol $P^n$ denotes the product $P \times \cdots \times P$ with $n$ factors $P$. If $X$ and $Y$ are disjoint, the {\em direct sum} $P + Q = ( X \cup Y, \leq_{P + Q} )$ is defined by
$$
\leq_{P+Q} \; := \; \leq_P \cup \leq_Q.
$$
The product $\ACh{n} \times P$, $n \in \myN$, is thus the direct sum of $n$ disjoint isomorphic copies of $P$.

If $P$ has a single maximal element, we call it the {\em top element} of $P$ and denote it by $\top_P$. Correspondingly, in the case of a single minimal element, we call it the {\em bottom element} and denote it by $\bot_P$.

If the reference poset is fixed or clear from the context, we skip the subscript ``$P$'' in notation.

\section{Counting down-sets} \label{sec_countingDownsets}

Let $M$ be a subset of $P$. For every $N \in \D(P \vert_M)$, we define
\begin{align*}
\D_{M,N}(P) & := \mysetdescr{ D \in \D(P) }{ D \cap M = N }.
\end{align*}
For $D \in \D(P)$, we have $D \cap M \in \D( P \vert_M )$. The set $\D(P)$ is thus the disjoint union of the sets $\D_{M,N}(P)$ with $N$ running through $\D(P \vert_M)$, and because of $ ( \darr_P N ) \cap M  = N$ for every $N \in \D(P \vert_M)$, the sets $\D_{M,N}(P)$ form a partition of $\D(P)$, thus
\begin{align} \label{dP_DNP_formel}
d(P) & = \sum_{N \in \D(P \vert_M)} \# \D_{M,N}(P).
\end{align}
As an example, for $\Ch{n} \times Q$, $n \in \myN$, we have with $M$ being the carrier of $\setx{0} \times Q$
$$
\D_{M,\setx{0} \times N}(\Ch{n} \times Q) \simeq \D \big( (\Ch{n-1}) \times Q \vert_N \big),
$$
for all $N \in \D( Q )$, hence \cite[Cor.\ 3.2]{Berman_Koehler_1976}
\begin{align} \nonumber
d( \Ch{n} \times Q ) & = \sum_{N \in \D(Q)} d \big( (\Ch{n-1}) \times Q \vert_N \big), \\ \label{d_PxChZwei}
\mytext{in particular,} \quad d( \Ch{2} \times Q) & = \sum_{N \in \D(Q)} \# \darr_{D(Q)} N.
\end{align}

Starting with the partition of $\D(P)$ formed by the sets $\D_{M,N}(P)$, $N \in \D( P \vert_M )$, an enumeration theorem for down-sets has been proven by the author and Ern\'{e} \cite[Theorem 3.1]{aCampo_Erne_2019} for $M$ being an {\em antichain} in $P$. In the following theorem, we generalize this result to {\em arbitrary subsets} $M$ of $P$. The advantage of the exotic sets $P - M \udarr_P N$ is explained later.

\begin{theorem} \label{theo_dP_formel}
Let $M$ be a subset of $P$. For every $N \in \D(P \vert_M)$, we define
\begin{align*}
M \udarr_P N & := \uarr_P \left( M \setminus N \right) \; \cup \; \darr_P N.
\end{align*}
The mapping
\begin{align*}
\phi_{M,N} : \D_{M,N}(P) & \rightarrow \D(P - M \udarr_P N), \\
D & \mapsto D \setminus \darr_P N,
\end{align*}
is an isomorphism for every $N \in \D(P \vert_M)$ with inverse
\begin{equation} \label{eq_phiN_invers}
\urbild{ \phi_{M,N} }(D') \; = \; D' \; \cup \; \darr_P N
\end{equation}
for all $D' \in \D(P - M \udarr_P N)$. In consequence,
\begin{align} \label{dP_formel}
d(P) & = \sum_{N \in \D(P \vert_M)} d( P - M \udarr_P N ).
\end{align}
\end{theorem}
\BP All arrows refer to $P$ and we skip the subscript ``$P$''.

Let $M$ be a fixed subset of $P$ and $N \in \D(P \vert_M)$ a fixed down-set of $P \vert_M$. Every down-set $D \in \D(P)$ satisfies the implications
$$
D \cap M = N \; \Rightarrow \; D \cap \uarr (M \setminus N) = \emptyset \; \Rightarrow \; ( D \setminus \darr N ) \cap ( M \udarr N) = \emptyset.
$$
Consequently, $D \in \D_{M,N}(P)$ implies $D \setminus \darr N \subseteq P - M \udarr N$. Furthermore, for every down-set $D \in \D(P)$, $x \in D$, and $y \in P - M \udarr N$  with $y \leq x$, the relation $y \notin \darr N$ trivially yields $y \in D \setminus \darr N$, and $\phi_{M,N}$ is a well-defined mapping. 

In order to see that $\phi_{M,N}$ is onto, let $E \in \D(P - M \udarr N)$, $D := E \cup \darr N$. Due to $E \in \D(P - M \udarr N)$, we have $E \cap \darr N = \emptyset$, hence $D\setminus \darr N = ( E \cup \darr N) \setminus \darr N = E$. It remains to show $D \in \D_{M,N}(P)$.

Assume $y \in D$ and $x \leq y$ with $x \notin D$. The relation $x \notin \darr N$ implies $y \notin \darr N$, hence $y \in E$, and $E \in \D( P - M \udarr N )$ in turn implies $x \in M \udarr N = \uarr (M \setminus N ) \cup \darr N$. But then $x \in \uarr( M \setminus N)$, thus $y \in \uarr( M \setminus N)$ in contradiction to $y \in E$. Therefore, $x \in D$, and $D$ is a down-set. Finally, due to $P - M \udarr N \subseteq P - M$,
$$
D \cap M = (E \cup \darr N) \cap M = ( \darr N ) \cap M = N,
$$
and $D \in \D_{M,N}(P)$ has been proven.

Now let $D \in \D_{M,N}(P)$. Due to $N \subseteq D \in \D(P)$, we have $\darr N \subseteq D$, hence $\phi_{M,N}(D) \cup \darr N = D$. The mapping $\phi_{M,N}$ is thus one-to-one with inverse \eqref{eq_phiN_invers}. Obviously, $\phi_{M,N}$ and its inverse are both order-preserving with respect to set-inclusion, and $\phi_{M,N}$ is an isomorphism. Now \eqref{dP_DNP_formel} yields \eqref{dP_formel}.

\EP

For the efficiency of formula \eqref{dP_formel}, the destructive power of subtracting $M \udarr N$ from $P$ is important. For $M$ being located somewhere in the mid of $P$, the set $\uarr (M \setminus N)$ will be large for small down-sets $N \in \D(P \vert_M)$, and $ \darr_P N$ will be large for large ones. We can thus expect that the down-set lattice $\D( P - M \udarr N)$ has moderate size for many $N \in \D(P \vert_M)$. In ideal case, we even can select $M$ in such a way that $P - M$ has a simple structure. Because of $P - M \udarr N \subseteq P - M$ for all $N \in \D(P \vert_M)$, this additionally facilitates the calculation of $d( P - M \udarr N )$. However, a too complicated set $M$ will blow up $\D( P \vert_M )$, causing a large number of summands.

In Section \ref{subsec_brrSechs}, we work with a poset $P = \Ch{2} \times Q$, and we need a result depending on this product structure only. With $Y$ being the carrier of $Q$, we define
\begin{align*}
M_0 & := \setx{0} \times Y, \quad P_0 := P\vert_{M_0}, \\
M_1 & := \setx{1} \times Y, \quad P_1 := P\vert_{M_1}.
\end{align*}
The mapping
\begin{align} \label{def_beta}
\begin{split}
\beta : P_0 &\rightarrow P_1,\\
(0,y) & \mapsto (1,y)
\end{split}
\end{align}
is an isomorphism. We have $x <_P \beta(x)$ for all $x \in M_0$, even
\begin{align} \label{x_leQ_y}
\forall \; x \in M_0, y \in M_1 \mytext{: }
\quad x <_P y \quad \Leftrightarrow \quad \beta(x) \leq_{P_1} y.
\end{align}

The result required in Section \ref{subsec_brrSechs} is

\begin{lemma} \label{lemma_QMzwodreiN}
For every $N \in \D(P_0)$,
\begin{align} \label{QMzwodreiN}
P - M_0 \udarr_P N = P_1 \vert_{\beta[ N ]}.
\end{align}
\end{lemma}
\BP We have $P - M_0 \udarr_P N \subseteq P - M_0 = M_1$, thus $P - M_0 \udarr_P N = M_1 \setminus \uarr_P ( M_0 \setminus N )$ because $N \subseteq M_0$ is a down-set in $P$, too. We show
$$
\urbild{\beta}\left( M_1 \setminus \uarr_P ( M_0 \setminus N ) \right) \; = \; N.
$$
For $y \in M_1 \setminus \uarr_P ( M_0 \setminus N )$, we must have $\urbild{\beta}(y) \notin M_0 \setminus N$ due to $\urbild{\beta}(y) <_P y$, hence $\urbild{\beta}(y) \in N$. On the other hand, for $y \in M_1 \cap \uarr_P ( M_0 \setminus N )$, there exists an $x \in M_0 \setminus N$ with $x <_P y$, hence $\beta(x) \leq_{P_1} y$ according to \eqref{x_leQ_y}. But then $x \leq_{P_0} \urbild{\beta}(y)$, thus $\urbild{\beta}(y) \in M_0 \setminus N$, because $M_0 \setminus N$ is an up-set in $P_0$.

\EP

\section{The Boolean lattice} \label{sec_BooleanLatt}

For $n \in \myN$, we define the {\em Boolean lattice $B(n)$ with $n$ atoms} as
$$
B(n) \; := \; \Ch{2}^n.
$$
The points of $B(n)$ are thus binary words with $n$ digits 0 or 1. For $x \in B(n)$, $x_1 \in \setx{0,1}$ is the first (leftmost) digit of $x$, and for every $0 \leq \ell \leq n$, we define the $\ell$\textsuperscript{th} {\em level} of $B(n)$ as
$$
L_\ell(n) := \mysetdescr{ x \in B(n) }{ x \mytext{ contains exactly } \ell \mytext{ ones} }.
$$
The elements of $L_1(n)$ are the {\em atoms} of $B(n)$ and the elements of $L_{n-1}(n)$ are the {\em co-atoms}. The only element of $L_0(n)$ is the bottom element $ \bot = 0 \cdots 0$ of $B(n)$, and the single element of $L_n(n)$ is the top element $ \top = 1 \cdots 1$.

Additionally, we need the following sub-posets of $B(n)$ for $n \geq2$:
\begin{align*}
\Brd(n) & := \cup_{i=2}^n L_i(n), \\
\Bru(n) & := \cup_{i=0}^{n-2} L_i(n), \\
\mytext {and} \quad \Brr(n) & := \cup_{i=2}^{n-2} L_i(n) \quad \mytext{ if } n \geq 3.
\end{align*}
$\Brd(n)$, $\Bru(n)$, and $\Brr(n)$ is thus the Boolean lattice with atoms and bottom element removed, with co-atoms and top element removed, and with atoms, co-atoms and both extrema removed, respectively. In order to unburden the notation, we additionally use the symbols
$$
b(n) := d( B(n) ), \quad
\brd(n) := d( \Brd(n) ), \quad
\bru(n) := d( \Bru(n) ), \quad
\brr(n) := d( \Brr(n) ).
$$

Even if the Dedekind numbers started their career as cardinalities of algebraic objects \cite{Dedekind_1897}, they turned out to be the number of down-sets (antichains, up-sets) of the Boolean lattices \cite[p.\ 61]{Birkhoff_1967}. The integer $b(n)$ is thus the $n$\textsuperscript{th} Dedekind number.

As pointed out by the author \cite{aCampo_2018}, all applicable algorithms for the calculation of Dedekind numbers described in literature \cite{Berman_Koehler_1976,aCampo_2018,Church_1940,
Church_1965,Fidytek_etal_2001,Lunnon_1971,
Markowsky_1989,Riviere_1968,Shmulevich_etal_1995,
Ward_1946,Wiedemann_1991,Yusun_2011} use the isomorphism $\D(B(n)) \simeq \H( B(\lambda), \D(B(n-\lambda)) )$ where $\H(P,P')$ is the poset formed by the order homomorphisms from $P$ to $P'$. Formula \eqref{d_PxChZwei} with $Q = B(n-1)$ results for $\lambda = 1$, and it is mentioned or used in many publications about Dedekind numbers \cite{Berman_Koehler_1976, aCampo_2018, Fidytek_etal_2001, Lunnon_1971, Markowsky_1989, Riviere_1968, Yusun_2011}. However, the standard algorithm used in the calculations \cite{Berman_Koehler_1976,Lunnon_1971,Markowsky_1989,
Wiedemann_1991,Yusun_2011} works with $\lambda = 2$ and it can be written as
\begin{align} \label{alg_L2_BK}
b(n) & = \sum_{ (D,E) \in \D(B(n-2))^2 } \# \darr_{\D(B(n-2))} ( D \cap E ) \cdot \# \uarr_{\D(B(n-2))} ( D \cup E ).
\end{align}
In enumerating the elements of $\H( B(2), \D(B(n-2)) )$, the summation runs over the images of the left and right corner of $B(2)$, and if they are fixed, the images of the top and bottom element of $B(2)$ are independently selected from the respective down-set and up-set in \eqref{alg_L2_BK}. Taking symmetry into account, the algorithm  requires $
\frac{b(n-2) \cdot ( b(n-2) + 1)}{2}
$
summands for the calculation of $b(n)$, hence 210 and 14196 summands for the calculation of $b(5)$ and $b(6)$, respectively. Furthermore, in a pre-calculation step, the integers $ \# \darr_{\D(B(n-2))} D$ and $\# \uarr_{\D(B(n-2))} D$ are calculated for every $D \in \D(B(n-2))$ and saved in a table.

Starting with Church \cite{Church_1940} in 1940, the numerous symmetries of the Boolean lattice have been used to reduce computational time. In particular, using a representation system of the non-isomorphic down-sets of $\D(B(n-2))$ speeds up the calculation considerably: For $b(7)$, Markowsky \cite{Markowsky_1989} reports a speed-up factor around 34. But of course, the pre-computational step to create a list of non-isomorphic elements of $\D( B(n-2) )$ is time consuming, and programming the main calculation becomes demanding \cite{Markowsky_1989,Wiedemann_1991}. More details about algorithms for the calculation of Dedekind numbers, including comparisons of their calculational effort, are found in literature \cite{aCampo_2018}.

\begin{figure}
\begin{center}
\includegraphics[trim = 70 260 80 70, clip]{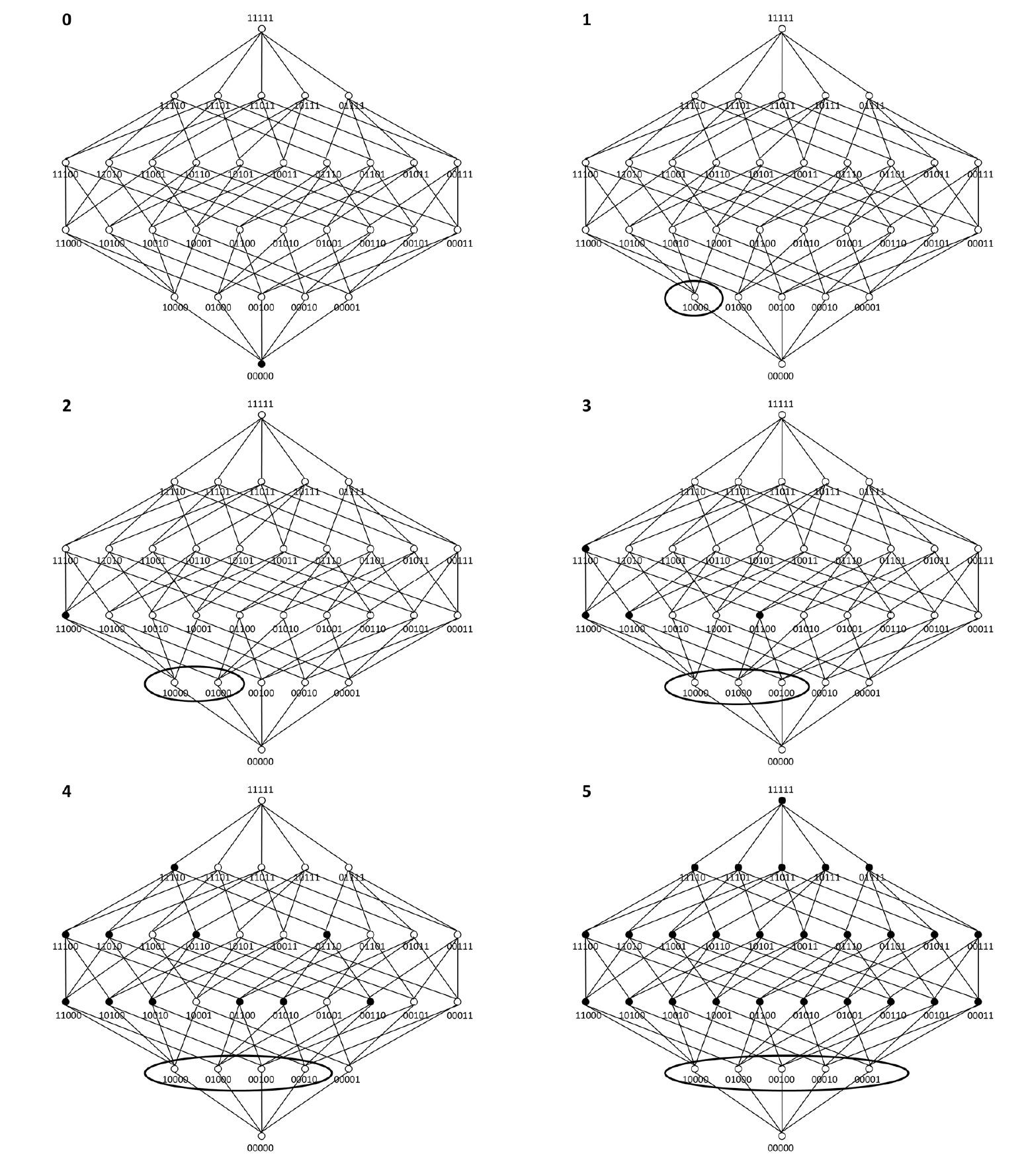}
\caption{\label{fig_Bn_minus_MN} The sets $B(5) - L_1(5) \udarr_{B(5)} N$ for different sizes of subsets $N \subseteq L_1(5)$. The sets $N$ are encircled and the points of the sets $B(5) - L_1(5) \udarr_{B(5)} N$ are shown as solid dots.}
\end{center}
\end{figure}

In the following theorem, we derive a formula for the Dedekind number $b(n)$ by applying Theorem \ref{theo_dP_formel} on $B(n)$ and $\Bru(n)$:

\begin{theorem} \label{theo_Dk}
For $n \geq 3$, we have
\begin{align} \label{dBn}
b(n) & = 2 + n + \sum_{k=2}^{n} \binom{n}{k} \cdot \brd(k).
\end{align}
Furthermore, $\brd(2) = 2$, and, for all $k \geq 3$,
\begin{align} \label{dBrd_n}
\brd(k) & = \brr(k) + 2 + \sum_{i=2}^{k-1} \binom{k}{i} \cdot \brd(i).
\end{align}
\end{theorem}
\BP Let $N \subseteq L_1(n)$, $k := \# N$. For $k = 0$, we have $\darr N = \emptyset$ and $\uarr( L_1(n) \setminus N) = B(n) \setminus \setx{ \bot_{B(n)} }$, hence $B(n) - L_1(n) \udarr N = \setx{ \bot_{B(n)} }$. For all other values of $k$, we have $\darr N = N \cup \setx{ \bot_{B(n)} }$, and the set $\uarr (L_1(n) \setminus N)$ contains exactly those points of $B(n)$ which have an 1 in a digit in which no element of $N$ has an 1. In consequence, the set $B(n) - L_1(n) \udarr_{B(n)} N$ contains exactly those points of $\Brd(n)$ which have all their 1s in the $k$ digits given by the 1s in the elements of $N$. All together (see Figure \ref{fig_Bn_minus_MN} for an illustration),
\begin{align} \label{Bn_subset}
B(n) - L_1(n) \udarr_{B(n)} N & \; \subseteq \; \bigcup_{i=2}^k L_i(n) \quad \;\; \; \mytext{for } 2 \leq k \leq n \\ \label{Bn_Beweis}
\mytext{and} \quad B(n) - L_1(n) \udarr_{B(n)} N & \; \simeq \; 
\begin{cases}
\setx{ \bot_{B(n)} }, & \mytext{if } k = 0; \\
\emptyset, & \mytext{if } k = 1; \\
\Brd( k ), & \mytext{if } 2 \leq k \leq n;
\end{cases}
\end{align}
and Theorem \ref{theo_dP_formel} yields \eqref{dBn}.

The equation $\brd(2) = 2$ is due to $\Brd(2) = \setx{ \top_{B(2)} }$. Furthermore, for all $k \geq 3$ and all $N \subseteq L_1(k)$,
$$
\Bru(k) - L_1(k) \udarr_{\Bru(k)} N \; \; = \; \; ( B(k) - L_1(k) \udarr_{B(k)} N ) \setminus ( L_{k-1}(k) \cup \setx{\top_{B(k)}} ),
$$
and \eqref{Bn_subset} and \eqref{Bn_Beweis} deliver with $i := \# N$
\begin{align} \label{Brun_AudN}
\Bru(k) - L_1(k) \udarr_{\Bru(k)} N & \simeq
\begin{cases}
\ACh{1}, & \mytext{if } i = 0; \\
\emptyset, & \mytext{if } i = 1; \\
\Brd( i ), & \mytext{if } 2 \leq i \leq k - 2; \\
\Brd( k-1 ) \setminus \setx{\top_{\Brd( k-1 )}}, & \mytext{if } i = k-1; \\
\Brr(k), & \mytext{if } i = k.
\end{cases}
\end{align}
Because $\Brd(k)$ and $\Bru(k)$ are dually isomorphic, Theorem \ref{theo_dP_formel} yields
\begin{align*}
\brd(k) & \; = \; \bru(k) \; = \;  \sum_{N \subseteq L_1(k)} d( \Bru(k) - L_1(k) \udarr_{\Bru(k)} N ),
\end{align*}
and \eqref{Brun_AudN} delivers \eqref{dBrd_n} because of $d \big( \Brd( k-1 ) \setminus \setx{\top_{\Brd( k-1 )}} \big) = \brd(k-1) - 1$.

\EP

\section{Application} \label{sec_Application}

\begin{table}
\begin{center}
\begin{tabular}{| c | r | r | r | }
\hline
$ n $ & $\brr(n)$ & $\brd(n)$ & $b(n)$ \\
\hline \hline
2 &    n.d. &        2 &       6 \\
3 &       1 &        9 &      20 \\
4 &      64 &      114 &     168 \\
5 &    6212 &     6894 &    7581 \\
6 & 7741776 &  7785062 & 7828354 \\
\hline
\end{tabular}
\caption{\label{table_dBn} The coefficients $\brr(n)$, $\brd(n)$, and the Dedekind numbers $b(n)$ for $2 \leq n \leq 6$.}
\end{center}
\end{table}

With \eqref{dBn} and \eqref{dBrd_n}, we have formulas allowing the recursive calculation of $b(n)$ and $\brd(n)$, provided that the values $\brr(k)$, $3 \leq k \leq n$, can been calculated. Due to $\Brr(3) = \emptyset$ and $\Brr(4) \simeq \ACh{6}$, we immediately get $\brr(3) = 1$ and $\brr(4) = 64$.


For the rest of the article, we focus on the application of Theorem \ref{theo_dP_formel} on the calculation of $\brr(5) = 6212$ and $\brr(6) = 7741776$ (Table \ref{table_dBn}). In the calculation of $\brr(5)$ and in the first brute force-approach to calculate $\brr(6)$, the set $M$ in Theorem \ref{theo_dP_formel} is an antichain. However, in the efficient second way to calculate $\brr(6)$ in Section \ref{subsec_brrSechs}, we use a subset $M \subset \Brr(6)$ for which $\Brr(6) \vert_M \simeq \Ch{2} \times \Brr(5)$ has a more complicated structure. As a by-product of the calculations, we get statistics about the sets $\D( \Brr(n) - M \udarr N )$ in the different methods. They are reported in the Tables \ref{table_aiFuenf}, \ref{table_ajcaFuenf}, and \ref{table_myijSechs}.

\subsection{The calculation of $\brr(5)$} \label{subsec_brrFuenf}

\begin{figure}
\begin{center}
\includegraphics[trim = 70 670 150 70, clip]{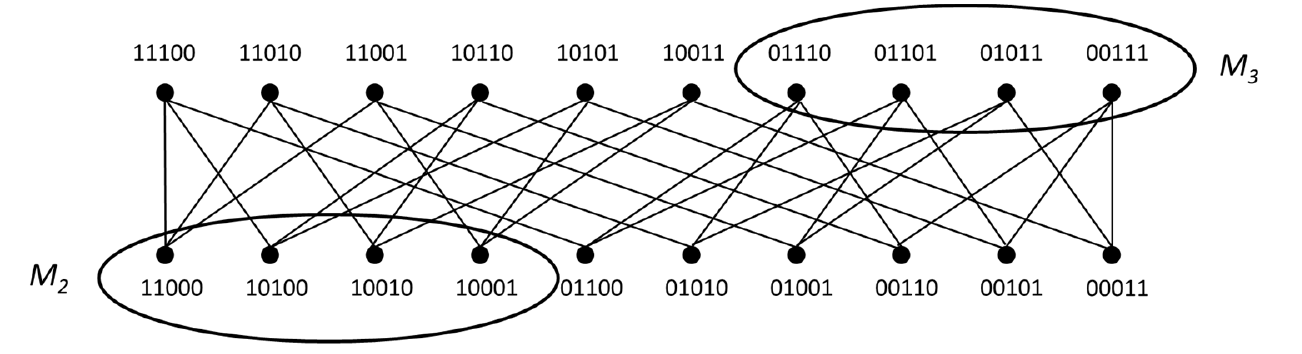}
\caption{\label{fig_L25L35} The poset $\Brr(5) = L_2(5) \cup L_3(5)$ and the sets $M_2$ and $M_3$ used in the second approach to calculate $\brr(5)$.}
\end{center}
\end{figure}

\begin{table}
\begin{center}
\begin{tabular}{| c || c | c | c | c | c | c | c | c | c | c | c |}
\hline
$i$ & 0 & 1 & 2 & 3 & 4 & 5 & 6 & 7 & 8 & 9 & 10 \\
\hline
$\nu_i$ & 388 & 290 & 195 & 70 & 40 & 30 & 0 & 10 & 0 & 0 & 1 \\
\hline
\end{tabular}
\caption{\label{table_aiFuenf} The coefficients $\nu_i$ used in \eqref{formel_brrFuenf} for the calculation of $\brr(5)$.}
\end{center}
\end{table}

The poset $\Brr(5) = L_2(5) \cup L_3(5)$ is shown in Figure \ref{fig_L25L35}. Each of the two level sets is an antichain with ten points. The most simple way to calculate $\brr(5)$ is to apply Theorem \ref{theo_dP_formel} with $M := L_3(5)$. The set $\D( \Brr(5) \vert_M )$ contains 1024 down-sets, and for each of them it is easy to calculate $\Brr(5) - M \udarr N = L_2(5) \setminus \darr N$ and $d( \Brr(5) - M \udarr N ) = 2^{\# ( L_2(5) \setminus \darr N )}$. With $\nu_i$ being the number of subsets $N \subseteq L_3(5)$ with $\# ( L_2(5) \setminus \darr N ) = i$, we have
\begin{equation} \label{formel_brrFuenf}
\brr(5) = \sum_{i=0}^{10} \nu_i \cdot 2^i.
\end{equation}
The coefficients $\nu_i$ are contained in Table \ref{table_aiFuenf}.

Even if the inspection of 1024 antichains $\Brr(5) - M \udarr N$ is not a big task, we want to reduce the effort. A natural idea is to use non-isomorphic down-sets in $M$, but this approach runs into problems because for $N \subseteq M$, the size of $L_2(5) \setminus \darr N$ is not uniquely determined by $\# N$ for $2 \leq \# N \leq 7$. In fact, we have to determine the 34 non-isomorphic down-sets of $\Brr(5)$ without isolated points, and for the pure calculation of $\brr(5)$, this effort does not pay. However, because we need these down-sets in Section \ref{subsec_brrSechs} for the calculation of $\brr(6)$, we have determined them, and in Formula \eqref{brrFümf_isomKl} in Section \ref{subsubsec_isom}, we manage the calculation of $\brr(5)$ by running a summation over 34 down-sets only.

For the calculation of $\brr(5)$ without determination of non-isomorphic down-sets, the following approach is efficient. We define
\begin{align*}
M_2 & := \mysetdescr{ x \in L_2(5) }{ x_1 = 1 }, \\
M_3 & := \mysetdescr{ x \in L_3(5) }{ x_1 = 0 },
\end{align*}
as illustrated in Figure \ref{fig_L25L35}. The set $M := M_2 \cup M_3$ is an eight-point antichain in $\Brr(5)$, and the set $\Brr(5) - M \simeq \ACh{6} \times \Ch{2}$ is the direct sum of six two-point-chains $\Ch{2}$. In consequence, for all $N \subseteq M$, the set $\Brr(5) - M \udarr N$ is the direct sum of at most six chains $\Ch{2}$ and at most six isolated points - it is thus easy to evaluate. Moreover, if $N_2, N_2'$ are subsets of $M_2$ with the same number of points, the sets 
\begin{align*}
& \mysetdescr{ \Brr(5) - M \udarr ( N_2 \cup N_3 ) }{ N_3 \subseteq M_3 } \\
\mytext{and} \quad
& \mysetdescr{ \Brr(5) - M \udarr ( N_2' \cup N_3 ) }{ N_3 \subseteq M_3 }
\end{align*}
contain for every class of non-isomorphic sub-posets of $\Brr(5) - M$ the same number of instances. It suffices thus to evaluate 80 posets $\Brr(5) - M \udarr ( N_2 \cup N_3 )$: Five sets $N_2 \subseteq M_2$ covering the possible cardinalities in combination with the sixteen subsets $N_3$ of $M_3$.

\begin{table}
\begin{center}
\begin{tabular}{| c || c | c | c | c | c | c || c | c | c || c | c | c | }
\hline
$ c $
& 0 & 0 & 0 & 0 & 0 & 0 & 1 & 1 & 1 & 2 & 3 & 6 \\ 
\hline
$ j \setminus a $
& 0 & 1 & 2 & 3 & 4 & 6 & 0 & 2 & 5 & 2 & 3 & 0 \\
\hline
\hline
 0 & 5 & 6 &   & 4 &   & 1 &   &   &   &   &   &   \\
 1 & 5 & 6 &   & 4 &   & 1 &   &   &   &   &   &   \\
 2 &   & 5 & 5 &   & 2 &   & 1 & 2 & 1 &   &   &   \\
 3 &   &   &   & 5 & 3 & 1 &   & 3 &   & 3 & 1 &   \\
 4 &   &   &   &   &   & 5 &   &   & 6 &   & 4 & 1 \\
\hline
\end{tabular}
\caption{\label{table_ajcaFuenf} The coefficients $\gamma_j(c,a)$ used in \eqref{formel_brrFuenf_Var} for the calculation of $\brr(5)$.}
\end{center}
\end{table}

For a fixed subset $N_2 \subseteq M_2$ of cardinality $j$, we define $\gamma_j(c, a)$ as the number of posets $N_3 \subseteq M_3$ with
$$
\Brr(5) - M \udarr ( N_2 \cup N_3 ) \; \simeq \; ( \ACh{c} \times\Ch{2} ) + \ACh{a}.
$$
($c$ is thus the number of 2-point chains in $\Brr(5) - M \udarr ( N_2 \cup N_3 )$ and $a$ is the number of isolated points.) Then
\begin{equation}  \label{formel_brrFuenf_Var}
\brr(5) = \sum_{c=0}^6 3^c \cdot \sum_{a=0}^6 2^a \cdot \sum_{j=0}^4 \binomk{4}{j} \cdot \gamma_j(c,a).
\end{equation}
The coefficients $\gamma_j(c,a)$ are shown in Table \ref{table_ajcaFuenf}.

\begin{figure}
\begin{center}
\includegraphics[trim = 70 460 120 70, clip]{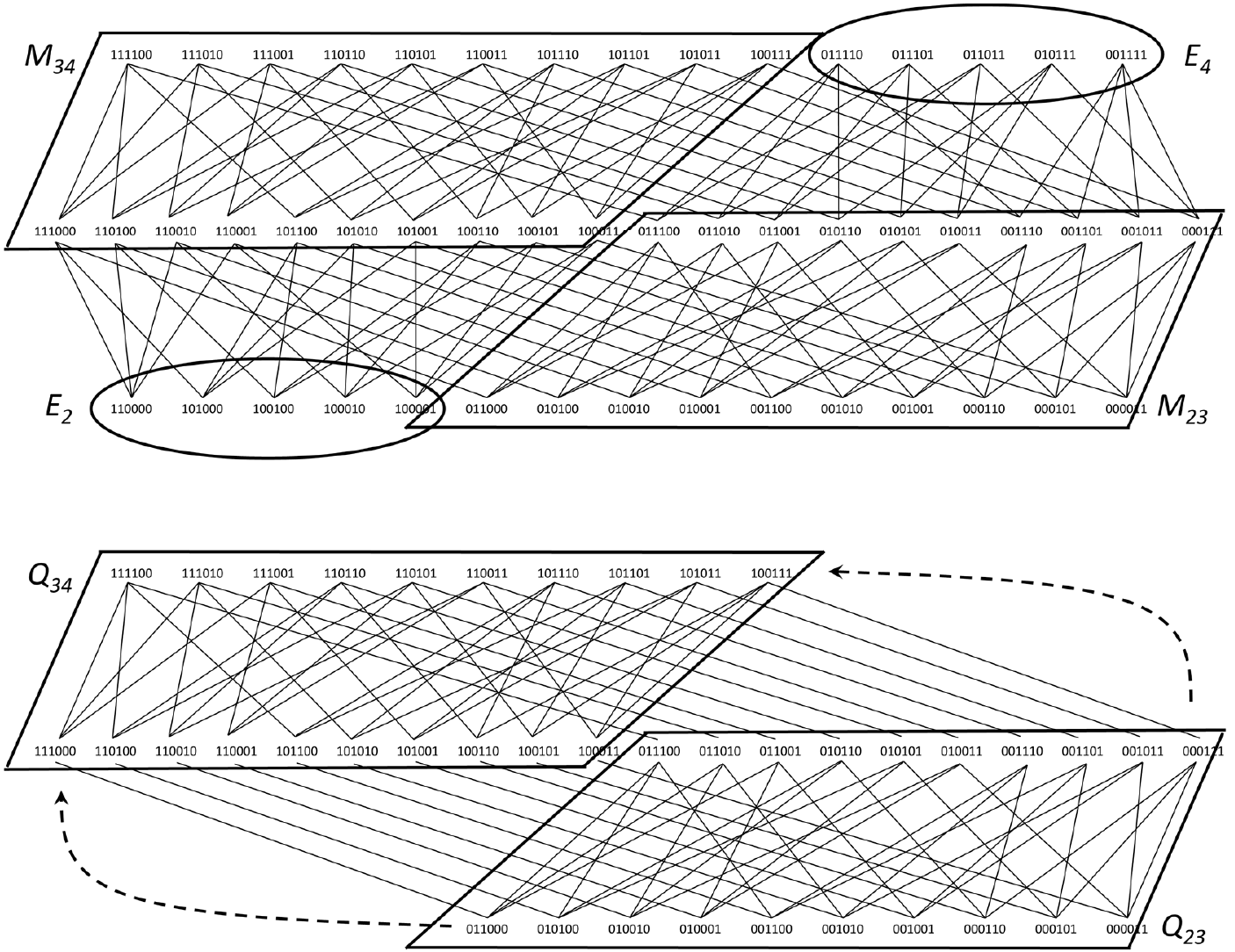}
\caption{\label{fig_L263646} Illustration of the structures used in the calculation of $\brr(6)$. In the upper part, the poset $\Brr(6) = L_2(6) \cup L_3(6) \cup L_4(6)$ is shown together with the sets $M_{23}$, $M_{34}$, $E_2$, and $E_4$. In the lower part, we have the poset $Q \simeq \Ch{2} \times \Brr(5)$ induced on $M_{23} \cup M_{34}$. The arrows indicate the isomorphism $\beta : Q_{23} \rightarrow Q_{34}$ from Lemma \ref{lemma_QMzwodreiN} for two points $x \in M_{23}$.}
\end{center}
\end{figure}

\subsection{The calculation of $\brr(6)$} \label{subsec_brrSechs}

\subsubsection{Approach} \label{subsubsec_approachSechs}

We come to the calculation of $\brr(6)$. Now the poset
$$
\Brr(6) = L_2(6) \cup L_3(6) \cup L_4(6)
$$
consists of three layers, as shown in Figure \ref{fig_L263646}. The mid-layer $L_3(6)$ is an antichain with 20 points, and even if it brings ordinary table calculation close to its limit, it is possible to calculate $\brr(6)$ by Formula \eqref{dP_formel} with $M := L_3(6)$ and $N$ running through all $2^{20}$ down-sets of $M$. (Again, the restriction to isomorphism classes is not possible without difficulties.)

Fortunately, the evaluation of $\Brr(6) - M \udarr N$ is simple. For $x \in L_2(6)$, $y \in L_4(6)$ with $x \leq y$, there exists an $m \in M$ with $x \leq m \leq y$, and therefore $x \notin \Brr(6) - M \udarr N$ or $y \notin \Brr(6) - M \udarr N$. The set $\Brr(6) - M \udarr N$ is thus always an antichain with ``lower'' part $P_2(N) := L_2(6) \cap ( \Brr(6) - M \udarr N )$ and ``upper'' part $P_4(N) := L_4(6) \cap ( \Brr(6) - M \udarr N )$. Denoting with $\mu(i,j)$ the number of subsets $N \subseteq M$ with $\# P_2(N) = i$ and $\# P_4(N) = j$, we have
\begin{equation} \label{formel_myijSechs}
\brr(6) = \sum_{i=0}^{15} \sum_{j=0}^{15} \mu(i,j) \cdot 2^{i+j}.
\end{equation}
The coefficients $\mu(i,j)$ are shown in Table \ref{table_myijSechs}.

\begin{table}
\begin{center}
{\scriptsize
\begin{tabular}{|r||r r r r|r r r r|r r r r|r r r r|}
\hline
$i \setminus j$&0&1&2&3&4&5&6&7&8&9&10&11&12&13&14&15 \\
\hline
\hline
0&165980&152265&86130&43385&17700&7569&2895&1350&420&160&90
&&20&&&1 \\
1&152265&103500&43080&16320&4410&1560&420&180&&15&&&&&&\\
2&86130&43080&13260&3660&585&180&60&&&&&&&&& \\
3&43385&16320&3660&800&&60&&&&&&&&&& \\
\hline
4&17700&4410&585&&&&&&&&&&&&& \\
5&7569&1560&180&60&&6&&&&&&&&&& \\
6&2895&420&60&&&&&&&&&&&&&\\
7&1350&180&&&&&&&&&&&&&&\\
\hline
8&420&&&&&&&&&&&&&&& \\
9&160&15&&&&&&&&&&&&&&\\
10&90&&&&&&&&&&&&&&& \\
11&&&&&&&&&&&&&&&& \\
\hline
12&20&&&&&&&&&&&&&&& \\
13&&&&&&&&&&&&&&&& \\
14&&&&&&&&&&&&&&&& \\
15&1&&&&&&&&&&&&&&& \\
\hline
\end{tabular}
}
\caption{\label{table_myijSechs} The coefficients $\mu(i,j)$ used in \eqref{formel_myijSechs} for the calculation of $\brr(6)$.}
\end{center}
\end{table}

We want to reduce the effort for the calculation of $\brr(6)$ by a smarter way to apply Theorem \ref{theo_dP_formel}. As indicated in Figure \ref{fig_L263646}, we define $M := M_{23} \cup M_{34}$ with 
\begin{align*}
M_{23} & := \mysetdescr{ x \in L_2(6) \cup L_3(6) }{ x_1 = 0 }, \\
M_{34} & := \mysetdescr{ x \in L_3(6) \cup L_4(6) }{ x_1 = 1 }
\end{align*}
and we want to calculate $\brr(6)$ by running $N$ in Theorem \ref{theo_dP_formel} through $\D( \Brr(6) \vert_M )$. In order to unburden the notation, we define $Q := \Brr(6) \vert_M$, $Q_{23} := Q \vert_{M_{23}} = \Brr(6) \vert_{M_{23}}$, and $Q_{34} := Q \vert_{M_{34}} = \Brr(6) \vert_{M_{34}}$.

At the first glance, the approach does not look promising because the posets $Q_{23}$ and $Q_{34}$ are both isomorphic to $\Brr(5)$, thus $Q \simeq \Ch{2} \times \Brr(5)$. We know that $\Brr(5)$ contains $6212 \approx 2^{12.6}$ down-sets, and we may argue that $Q$ has clearly more than the $2^{20}$ down-sets evaluated in the brute force-approach. (Indeed, the number is $3933651 \approx 2^{21.9}$.) However, we can reduce the effort if we are willing to invest into the analysis of $\D(Q)$. We show in Lemma \ref{lemma_brrSechs} that we  can calculate $\brr(6)$ by a summation running over the 6212 down-sets contained in $\D(Q_{23})$, and focusing on non-isomorphic down-sets, we manage the calculation in Section \ref{subsubsec_isom} with a summation over 245 down-sets.

We need some additional definitions. Setting as in Figure \ref{fig_L263646}
\begin{align*}
E_2 & := \mysetdescr{ x \in L_2(6) }{ x_1 = 1 }, \\
E_4 & := \mysetdescr{ x \in L_4(6) }{ x_1 = 0 },
\end{align*}
$E_2$ and $E_4$ are disjoint five-point antichains in $\Brr(6)$ with $\Brr(6) - M = E_2 + E_4$ being a ten-point antichain in $\Brr(6)$. We define
\begin{align*}
S & := Q \vert_{E_2 \cup M_{34}}, \\
T & := Q \vert_{E_4 \cup M_{23}}.
\end{align*}
$S$ is thus the left part of the poset $\Brr(6)$ in Figure \ref{fig_L263646} whereas $T$ is its right part.

The points of $M_{23}$ and $E_2$ are pairwise incomparable in $\Brr(6)$, just as the points of $M_{34}$ and $E_4$. Therefore, for all $N \in \D(Q)$,
$$
\Brr(6) - M \udarr_{\Brr(6)} N
\;\; = \;\; 
E_2 \setminus \darr_S (M_{34} \cap N)
\;\; \cup \;\;
E_4 \setminus \uarr_T (M_{23} \setminus N),
$$
hence
\begin{align}  \label{zweisNtN}
d( \Brr(6) - M \udarr_{\Brr(6)} N ) & = 2^{s(N) + t(N)} \\ \label{formel_sN}
\mytext{with} \quad s(N) & := \# \big( E_2 \setminus \darr_S (M_{34} \cap N) \big) \\ \label{formel_tN}
\mytext{and} \quad t(N) & := \# \big( E_4 \setminus \uarr_T (M_{23} \setminus N) \big). \\ \label{tN_tN_cap_Mzd}
\mytext{In particular,} \quad t(N) & \; = t(N \cap M_{23} ).
\end{align}

We have $Q \simeq \Ch{2} \times \Brr(5)$ and we can therefore apply Lemma \ref{lemma_QMzwodreiN} on $Q$. The isomorphism $\beta : Q_{23} \rightarrow Q_{34}$ used in the lemma and defined in \eqref{def_beta} is given by switching the first digit of $x \in M_{23}$ from 0 to 1 (cf.\ Figure \ref{fig_L263646}).

\begin{lemma} \label{lemma_brrSechs}
\begin{align} \label{summe_brrSechs}
\brr(6) & = \sum_{N \in \D(Q_{23})} 2^{t(N)} \cdot \sigma(N) \\ \label{def_sigmaSumme}
\mytext{with} \quad
\sigma(N) & := \sum_{N' \in \darr_{\D(Q_{23})} N} 2^{\# ( E_2 \setminus \darr_S \beta[N'])}
\quad \mytext{for all } N \in \D( Q_{23} ).
\end{align}
\end{lemma}
\BP In the beginning of Section \ref{sec_countingDownsets}, we have seen that the sets
$$
\D_{M_{23},N}(Q) = \mysetdescr{ D \in \D(Q) }{ D \cap M_{23} = N }, \quad N \in \D( Q_{23} ),
$$
form a partition of $\D(Q)$, and according to Theorem \ref{theo_dP_formel}, the mapping
\begin{align*}
\phi_{M_{23},N} : \D_{M_{23},N}(Q) & \rightarrow \D( Q - M_{23} \udarr_Q N ), \\
D & \mapsto D \setminus \darr_Q N
\end{align*}
is an isomorphism for all $N \in \D(Q_{23})$. Therefore, 
\begin{align*}
\brr(6)
& \streleq{dP_formel}
\sum_{N \in \D(Q)} d( \Brr(6) - M \udarr_{\Brr(6)} N ) \\
& \streleq{zweisNtN} 
\sum_{N \in \D(Q)} 2^{s(N) + t(N)} \\
& =
\sum_{N \in \D(Q_{23})} \; \sum_{D \in \D_{M_{23},N}(Q)} 2^{s(D) + t(D)} \\
& \streleq{tN_tN_cap_Mzd}
\sum_{N \in \D(Q_{23})} 2^{t(N)} \sum_{D \in \D_{M_{23},N}(Q)} 2^{s(D)} \\
& =
\sum_{N \in \D(Q_{23})} 2^{t(N)} \sum_{D' \in \D(Q - M_{23} \udarr_Q N )} 2^{s \left( \urbild{\phi_{M_{23},N}}(D') \right) } \\
& \streleq{QMzwodreiN} 
\sum_{N \in \D(Q_{23})} 2^{t(N)} \sum_{D' \in \D(Q_{34} \vert_{\beta[N]})} 2^{s \left( \urbild{\phi_{M_{23},N}}(D') \right) }.
\end{align*}
For $D' \in \D(Q_{34} \vert_{\beta[N]})$, $N \in \D(Q_{23})$,
\begin{align*}
M_{34} \cap \urbild{\phi_{M_{23},N}}(D')
& \; \; \streleq{eq_phiN_invers} \; \; 
M_{34} \cap \big( D' \; \cup \; \darr_Q N \big)
\; \; = \; \; D', \\
\mytext{thus} \quad 
s \left( \urbild{\phi_{M_{23},N}}(D') \right)  
& \; \; \streleq{formel_sN} \; \; 
\# ( E_2 \setminus \darr_S D' ),
\end{align*}
and because $\beta : Q_{23} \rightarrow Q_{34}$ is an isomorphism,
$$
\sum_{D' \in \D(Q_{34} \vert_{\beta[N]})} 2^{\# ( E_2 \setminus \darr_S D')}
\; \; = \; \;
\sum_{N' \in \D(Q_{23} \vert_N )} 2^{\# ( E_2 \setminus \darr_S \beta[N'] )}.
$$
Finally, $\D(Q_{23} \vert_N ) = \darr_{\D(Q_{23})} N$ for every $N \in \D(Q_{23})$.

\EP

In Formula \eqref{summe_brrSechs}, the summation runs over 6212 summands only. However, the calculation of the coefficients $\sigma(N)$ by means of \eqref{def_sigmaSumme} requires much effort. In the following section, we reduce the number of summands by taking isomorphism into account, and we derive a handy formula for $\sigma(N)$. As already announced, we finally manage the calculation of $\brr(6)$ with the quick calculation of $\sigma(N)$ for 245 down-sets $N \in \D(Q_{23})$.

\subsubsection{Exploiting isomorphism} \label{subsubsec_isom}

For a subset $Y \subseteq M_{23}$, we call the points in $Y \cap L_2(6)$ the {\em lower points} of $Y$ and the points in $Y \cap L_3(6)$ the {\em upper points} of $Y$.

\begin{figure}
\begin{center}
\includegraphics[trim = 75 340 70 70, clip]{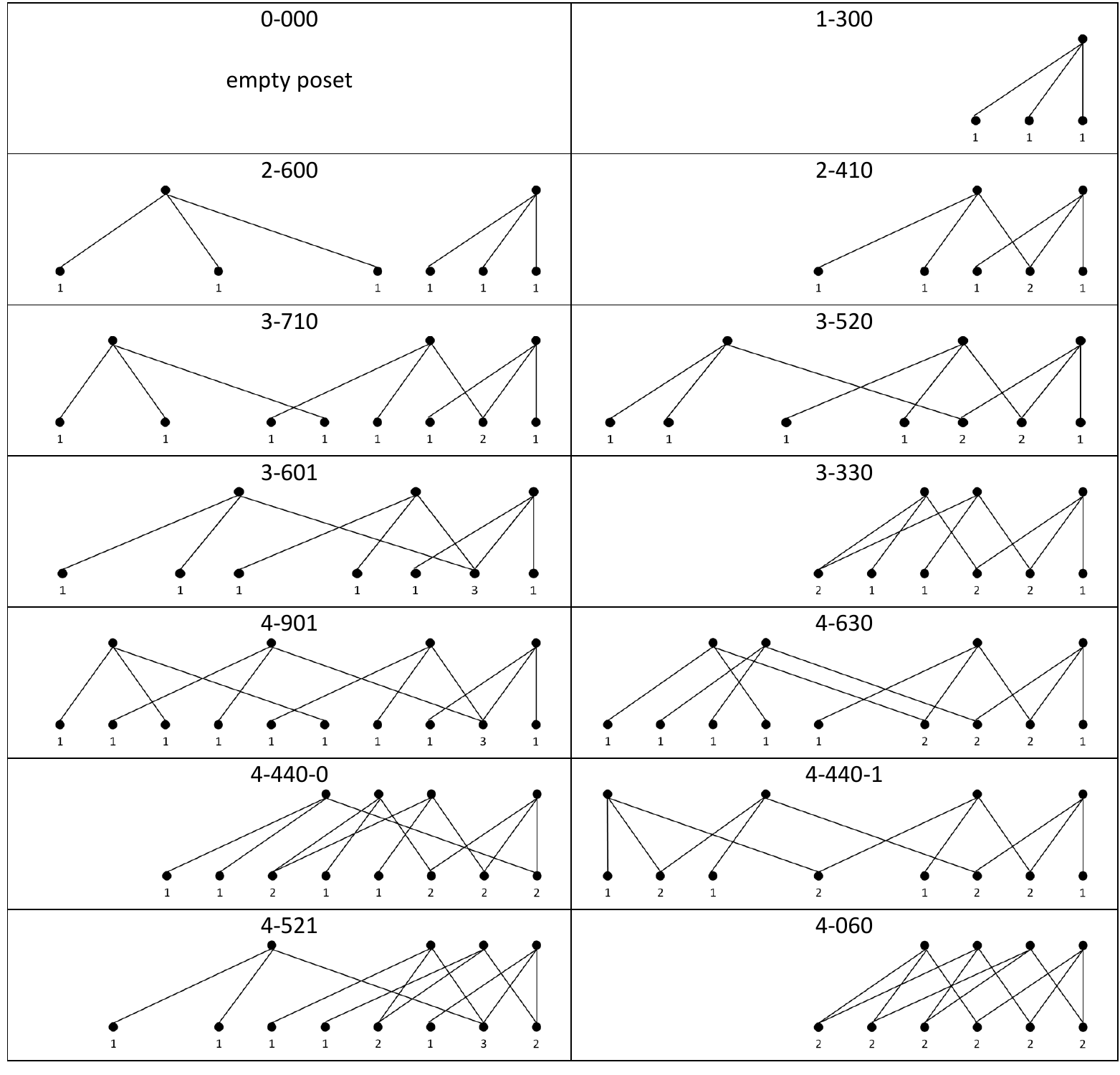}
\caption{\label{fig_IsomKl_01234} The non-isomorphic down-sets of $Q_{23} \simeq \Brr(5)$ with up to four upper points and no isolated points. The integers under the lower points indicate the number of their upper covers.}
\end{center}
\end{figure}

\begin{figure}
\begin{center}
\includegraphics[trim = 75 570 70 70, clip]{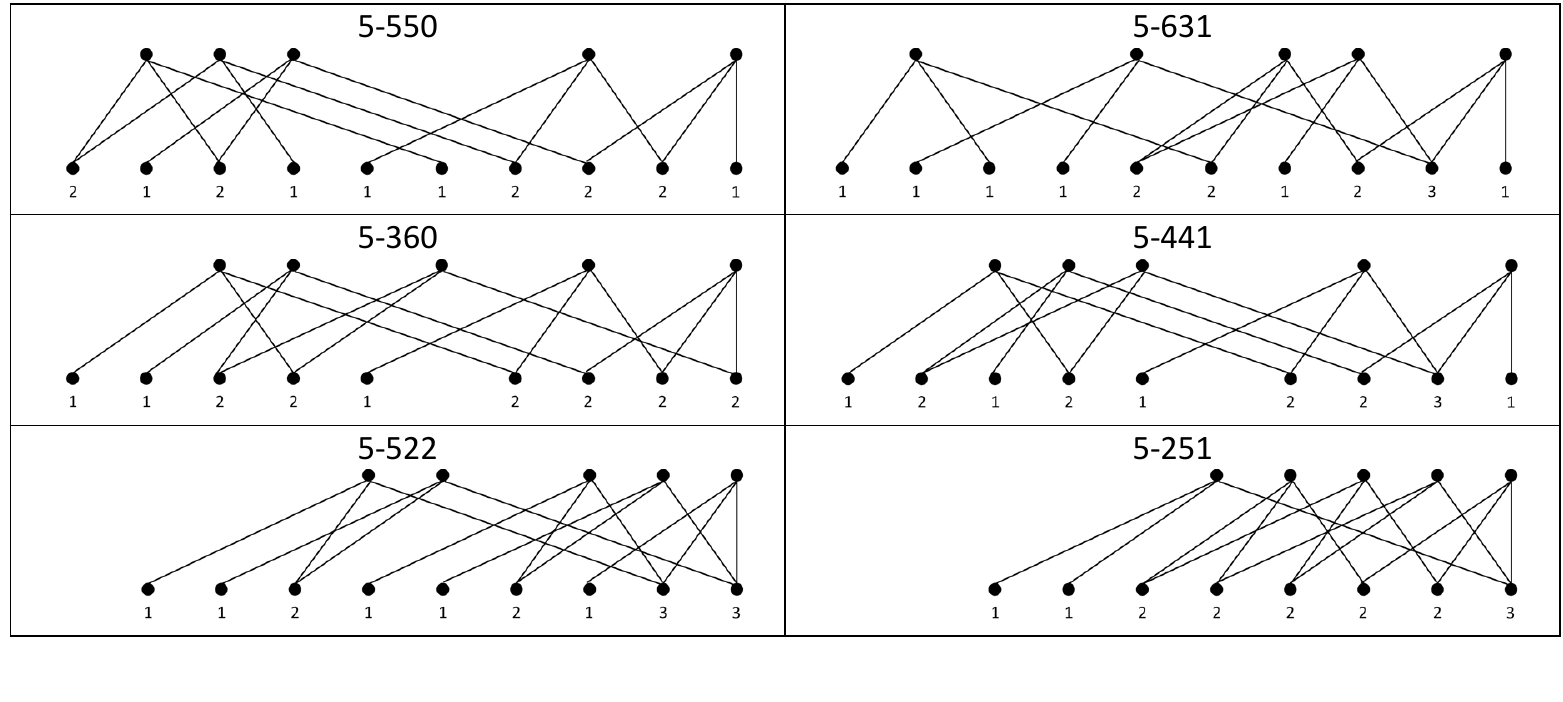}
\caption{\label{fig_IsomKl_5} The non-isomorphic down-sets of $Q_{23} \simeq \Brr(5)$ with five upper points and no isolated points. The integers under the lower points indicate the number of their upper covers.}
\end{center}
\end{figure}

\begin{figure}
\begin{center}
\includegraphics[trim = 75 340 70 70, clip]{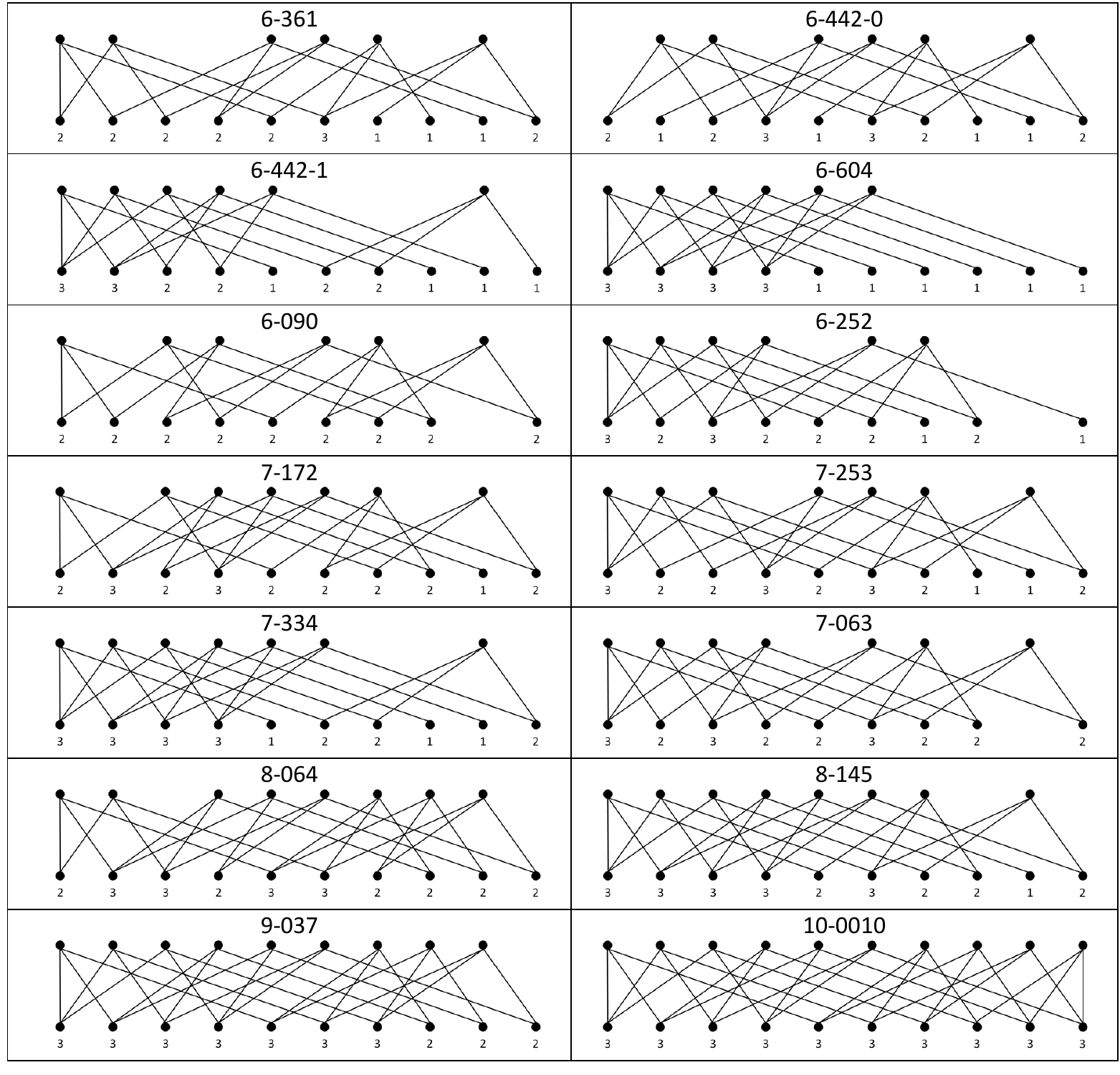}
\caption{\label{fig_IsomKl_678910} The non-isomorphic down-sets of $Q_{23} \simeq \Brr(5)$ with at least six upper points and no isolated points. The integers under the lower points indicate the number of their upper covers.}
\end{center}
\end{figure}

Let $\R$ denote a representation system of the non-isomorphic down-sets of $Q_{23} \simeq \Brr(5)$ and let $\R_0 \subset \R$ be the set of posets contained in $\R$ without isolated points. The thirty-four elements of $\R_0$ are shown in the Figures \ref{fig_IsomKl_01234}, \ref{fig_IsomKl_5}, and \ref{fig_IsomKl_678910}. Thirty of them are uniquely described by a four-integer code
$$
u\mytext{-}c_1c_2c_3,
$$
in which, for $R \in \R_0$, $u$ is the number of upper points in $R$ and $c_j$ is the number of lower points in $R$ covered by $j$ upper points. For the remaining four elements of $\R_0$, we make the code unique with an additional digit:
\begin{itemize}
\item 4-440-0: The down-set $R$ is of type 4-440 and does not contain an 8-crown.
\item 4-440-1: The down-set $R$ is of type 4-440 and contains an 8-crown.
\item 6-442-0: The down-set $R$ is of type 6-442 and for each upper point $x \in R$, there exists a lower point $y \in ( \darr x ) \setminus \setx{x}$ which is covered by three upper points.
\item 6-442-1: The down-set $R$ is of type 6-442 and there exists an upper point $x \in R$ for which no lower point $y \in ( \darr x ) \setminus \setx{x}$ is covered by three upper points.
\end{itemize}
Table \ref{table_TypeParas} below contains parameters for every $R \in \R_0$. The integer $\iota(R)$ is the number of isomorphic copies of $R$ in $Q_{23}$, and with $\Delta(R) := ( L_2(6) \cap M_{23} ) \setminus R$, the integer $\delta(R) := \# \Delta(R)$ is the number of lower points in $M_{23}$ not belonging to $R$. (For $R$ of type $u$-$c_1c_2c_3$, $\delta(R) = 10-c_1-c_2-c_3$.) The parameter $t(R)$ has already been introduced in \eqref{formel_tN}. Additionally, $\sigma(R)$, $\# \darr_{Q_{23}} R$ and the result of the inner summation in Formula \eqref{formel_brrSechs} below are shown. We have
\begin{align} \label{brrFümf_isomKl}
\brr(5) & = \sum_{R \in \R_0} \iota(R) \cdot 2^{\delta(R)} \\ \nonumber
\mytext{and} \quad \# \R & = \sum_{R \in \R_0}(  1 + \delta(R) ) \; = \; 91.
\end{align}

Even if the summation in \eqref{brrFümf_isomKl} is running over $\R_0$, it is in fact an application of Theorem \ref{theo_dP_formel} on $\Brr(5)$ with $M := L_3(5)$. Each $R \in \R_0$ is uniquely determined by its upper points $R \cap L_3(5)$, and we have $d( \Brr(5) - L_3(5) \udarr (R \cap L_3(5) ) ) = 2^{\delta( R )}$. Selecting an appropriate partition of the power set of $L_3(5)$, we can rewrite \eqref{brrFümf_isomKl} in such a way that it is fully in line with \eqref{dP_formel}.

With 91 elements, the representation system $\R$ is small, but it is of restricted value for the calculation of $\brr(6)$, because for $N \in \D(Q_{23})$, the cardinality of $E_2 \setminus \darr_S \beta[N]$ is not uniquely determined by the isomorphism type of $N$. However, according to Theorem \ref{theo_dP_formel} and Lemma \ref{lemma_brrSechs},
\begin{align} \label{formel_brrSechs}
\brr(6) & \; = \; \sum_{R \in \R_0} \iota(R) \cdot 2^{t(R)} \cdot \sum_{A \subseteq \Delta(R) } \sigma(A + R).
\end{align}

In this formula, the two nested sums create $\sum_{R \in \R_0} 2^{\delta(R)} = 1269$ down-sets of the form $A + R$, which is not much, but the calculation of the coefficients $\sigma(A + R)$ with the defining formula \eqref{def_sigmaSumme} requires all together 208099 evaluations. We have to develop a formula for the quick calculation of $\sigma(A + R)$. As short-cut, we define
$$
e(Y) \; := \; \# \big( E_2 \setminus \darr_S \beta[Y] \big) \quad \mytext{for all } Y \subseteq M_{23}.
$$

We start with calculating two tables with 1024 entries each: For each set $Y \subseteq L_2(6) \cap M_{23}$ of lower points of $M_{23}$, we set
\begin{align*}
\fT_0(Y) & := 2^{e(Y)}, \\
\fT_1(Y) & := \sum_{Z \subseteq Y} \fT_0(Z).
\end{align*}

Of the 5188 down-sets $D \in \D(Q_{23})$ with upper points, only 491 have $e(D) > 0$, and after removing their isolated points, all are of type 1-300, 2-410, 3-330, or 4-060. 

Now let $R \in R_0$ and $A \subseteq \Delta(R)$. We want to determine $\sigma(R')$ with $R' := A + R$. The set $A' := R' \cap L_2(6)$ is the set of lower points of $R'$, and $U := R \cap L_3(6)$ is the set of upper points of $R'$.

There are $2^{\# A'}$ down-sets $D \in \darr_{\D(Q_{23})} R'$ without upper points, and their total contribution to $\sigma(R')$ is
$$
\sum_{
\stackrel{D \in \darr_{\D(Q_{23})} R'}{D \cap L_3(6) = \emptyset}
} 2^{e(D)} \; \; = \; \; \fT_1(A').
$$

For each of the types $\tau \in \setx{ \mytext{1-300, 2-410, 3-330, 4-060} }$, let $\rho_\tau$ be the number of down-sets $D \in \darr_{\D(Q_{23})} R'$ with $s(D) > 0$ for which $D$ with isolated points removed is of type $\tau$. These numbers can be determined as follows:
\begin{itemize}
\item For type 1-300, we step through the points $u \in U$. For each $u \in U$, there exist two disjoint three-point sets $G_1(u), G_2(u) \subseteq \Delta( \darr_{Q_{23}} u )$ with the following property: For all $Y \subseteq \Delta( \darr_{Q_{23}} u )$,
\begin{align*}
e(Y \cup \darr_{Q_{23}} u) & = 
\begin{cases}
2, & \mytext{if } Y = \emptyset; \\
1, & \mytext{if } Y \not= \emptyset \mytext{ and} Y \subseteq G_1(u) \mytext{ or} Y \subseteq G_2(u); \\
0, & \mytext{otherwise.}
\end{cases}
\end{align*}
Therefore, each $u \in U$ contributes to $\rho_{\mytext{1-300}}$ with $2^{\# G_1(u) \cap A'} + 2^{\# G_2(u) \cap A'} - 1$ items.
\item For each two-element subset $V \subseteq U$, we check the down-set $\darr_{Q_{23}} V$. It is of type 2-410 iff it contains exactly five lower points. In this case,  there exists a single point $g(V) \in \Delta( \darr_{Q_{23}} V ) $ with the following property: For all $Y \subseteq \Delta( \darr_{Q_{23}} V )$,
\begin{align*}
e(Y \cup \darr_{Q_{23}} V) & = 
\begin{cases}
1, & \mytext{if } Y \subseteq \setx{ g(V) }; \\
0, & \mytext{otherwise.}
\end{cases}
\end{align*}
Each two-element subset $V \subseteq U$ with $\darr_{Q_{23}} V$ being of type 2-410 contributes thus to $\rho_{\mytext{2-410}}$ with two items if $g(V) \in A'$ and with a single item otherwise.

\item For the types 3-330 and 4-060, we step through the subsets $V \subseteq U$ with three and four elements, respectively, and we check if the resulting down-set $\darr_{Q_{23}} V$ is of the respective type. That is easily done, because $\darr_{Q_{23}} V$ is of type 3-330 or 4-060, respectively, iff it contains exactly six lower points. For all $Y \subseteq \Delta( \darr_{Q_{23}} V )$, we have
\begin{align*}
e(Y \cup \darr_{Q_{23}} V) & = 
\begin{cases}
1, & \mytext{if } Y = \emptyset; \\
0, & \mytext{otherwise.}
\end{cases}
\end{align*}
Each $V \subseteq U$ with three and four elements and $\darr_{Q_{23}} V$ being of type 3-330 or 4-060 contributes thus with a single item to $\rho_{\mytext{3-330}}$ and $\rho_{\mytext{4-060}}$, respectively.
\end{itemize}
Now we can calculate $\sigma(R')$ as follows:
\begin{itemize}
\item The table entry $\fT_1(A')$ is the contribution of the $2^{\# A'}$ elements of $\darr_{Q_{23}} R' $ without upper points.
\item There are $\rho_{\mytext{1-300}} + \rho_{\mytext{2-410}} + \rho_{\mytext{3-330}} + \rho_{\mytext{4-060}}$ elements $D \in \darr_{Q_{23}} R' $ with upper points and $e(D) > 0$. All together, they contribute
$$
2 \cdot \left( \# U + \rho_{\mytext{1-300}} + \rho_{\mytext{2-410}} + \rho_{\mytext{3-330}} + \rho_{\mytext{4-060}} \right).
$$
(The term ``$\# U \; + \; $'' is due to $e( \darr_{Q_{23}} u ) = 2$ for each $u \in U$.)
\item Of the remaining
$$
\# \left( \darr_{Q_{23}} R' \right) - 2^{\# A'} - \rho_{\mytext{1-300}} - \rho_{\mytext{2-410}} - \rho_{\mytext{3-330}} - \rho_{\mytext{4-060}}
$$
elements of $\darr_{Q_{23}} R'$, each contributes 1.
\item Finally, with $\# \darr_{Q_{23}} R$ from Table \ref{table_TypeParas},
$$
\# \darr_{Q_{23}} R' \quad = \quad 2^{\# A} \cdot \# \darr_{Q_{23}} R.
$$ 
\end{itemize}
Putting all together yields
\begin{align} \label{formel_sigma}
\begin{split}
\sigma(R')
\quad = & \quad \; \; \; \; \left( 2^{\# A} \cdot \# \darr_{Q_{23}} R \right) - 2^{\# A'} + \fT_1(A') \\
& \quad + 2 \cdot \# U + \rho_{\mytext{1-300}} + \rho_{\mytext{2-410}} + \rho_{\mytext{3-330}} + \rho_{\mytext{4-060}}.
\end{split}
\end{align}

In calculating $\brr(6)$ with Formula \eqref{formel_brrSechs}, we can skip the generation and evaluation of down-sets without upper points: their contribution to $\brr(6)$ is the sum over all entries in the table $\fT_1$. It remains to generate and evaluate the down-sets $A + R$ with $R \in \R_0$ having upper points and $A \subseteq \Delta(R)$. If we regard the creation of the table $\fT_1$ as belonging to pre-calculation, we manage the main calculation of $\brr(6)$ with evaluating $1269 - 1024 = 245$ down-sets.

For the calculation of $\sigma(R')$ with Formula \eqref{formel_sigma}, parts of the analyses we have to do for the determination of the numbers $\rho_\tau$ depend on $R$ only. For each $R \in \R_0$ in the outer summation in \eqref{formel_brrSechs}, they have to be done only once for the calculation of all coefficients $\sigma(A + R)$ with $A \subseteq \Delta(R)$.

Using the information about types, isolated points and values of $e(A + R)$ collected in this section, it is possible to calculate $\brr(6)$ efficiently with a summation running over the 91 isomorphism classes $R' \in \R$ of down-sets of $Q_{23}$, even with a summation running over the 34 elements of $\R_0$ only. However, in our opinion, the calculation of the value of the summands becomes complicated and violates thus the second and third criterion of simplicity set up in the introduction. Therefore, we did not calculate $\brr(6)$ in this way.

\begin{table}
\begin{center}
{\scriptsize
\begin{tabular}{| l || r | r | r | r | r | r |}
\hline
type of $R$ & $\iota(R)$ & $\delta(R)$ & $t(R)$ & $\sigma(R) $ & $\# \darr_{\D(Q_{23})} R$ & $\sum_{A \subseteq \Delta(R)} \sigma(A + R)$\\
\hline
\hline
0-000 & 1 & 10 & 0 & 32 & 1 & 173433\\
\hline
1-300 & 10 & 7 & 0 & 76 & 9 & 42075\\
\hline
2-600 & 15 & 4 & 0 & 221 & 81 & 10821\\
2-410 & 30 & 5 & 0 & 166 & 41 & 17711\\
\hline
3-710 & 30 & 2 & 0 & 644 & 369 & 4791\\
3-520 & 60 & 3 & 0 & 387 & 187 & 7621\\
3-601 & 10 & 3 & 0 & 403 & 189 & 7738\\
3-330 & 20 & 4 & 0 & 294 & 95 & 12481\\
\hline
4-901 & 10 & 0 & 0 & 2201 & 1701 & 2201\\
4-630 & 60 & 1 & 0 & 1227 & 853 & 3433\\
4-440-0 & 60 & 2 & 0 & 728 & 434 & 5462\\
4-440-1 & 15 & 2 & 0 & 697 & 433 & 5413\\
4-521 & 60 & 2 & 0 & 736 & 439 & 5519\\
4-060 & 5 & 4 & 1 & 332 & 113 & 14297\\
\hline
5-550 & 12 & 0 & 0 & 2496 & 1975 & 2496\\
5-631 & 60 & 0 & 0 & 2530 & 2006 & 2530\\
5-360 & 60 & 1 & 0 & 1400 & 1007 & 3938\\
5-441 & 60 & 1 & 0 & 1423 & 1022 & 3994\\
5-522 & 30 & 1 & 0 & 1437 & 1035 & 4036\\
5-251 & 30 & 2 & 1 & 842 & 524 & 6378\\
\hline
6-361 & 60 & 0 & 0 & 2925 & 2377 & 2925\\
6-442-0 & 15 & 0 & 1 & 2984 & 2431 & 2984\\
6-442-1 & 60 & 0 & 0 & 2967 & 2416 & 2967\\
6-604 & 5 & 0 & 0 & 3045 & 2489 & 3045\\
6-090 & 10 & 1 & 0 & 1607 & 1195 & 4545\\
6-252 & 60 & 1 & 1 & 1666 & 1241 & 4704\\
\hline
7-172 & 30 & 0 & 0 & 3456 & 2881 & 3456\\
7-253 & 60 & 0 & 1 & 3529 & 2949 & 3529\\
7-334 & 20 & 0 & 1 & 3584 & 3001 & 3584\\
7-063 & 10 & 1 & 2 & 1968 & 1519 & 5591\\
\hline
8-064 & 15 & 0 & 1 & 4214 & 3607 & 4214\\
8-145 & 30 & 0 & 2 & 4310 & 3698 & 4310\\
\hline
9-037 & 10 & 0 & 3 & 5337 & 4693 & 5337\\
\hline
10-0010 & 1 & 0 & 5 & 6893 & 6212 & 6893\\
\hline
\end{tabular}
\caption{\label{table_TypeParas} Parameters of the non-isomorphic down-sets $R \in \R_0$ without isolated points. Explanation in text.}
}
\end{center}
\end{table}

(Concerned with sequence \seqnum{A000372}.)


\begin{thebibliography}{xx}

\bibitem{Berman_Koehler_1976} J. Berman and P. K\"{o}hler, Cardinalities of finite distributive lattices, {\em Mitteilungen aus dem mathem.\ Seminar Gie{\ss}en} {\bf 121} (1976), 103--124.

\bibitem{Birkhoff_1967} G. Birkhoff, {\em Lattice Theory}, Proc.\ Amer.\ Math.\ Soc.\ Coll.\ Publ.\ {\bf 25}, 3\textsuperscript{rd} ed., 1967.

\bibitem{aCampo_2018} F. a Campo, Relations between powers of Dedekind numbers and exponential sums related to them, {\em J.\ Integer Seq.} {\bf 21} (2018), article 18.4.4, \url{https://cs.uwaterloo.ca/journals/JIS/VOL21/Campo/campo3.pdf} and \url{https://cs.uwaterloo.ca/journals/JIS/VOL21/Campo/campo3-corrigendum.pdf}.

\bibitem{aCampo_2018_Framework} F. a Campo, A framework for the systematic determination of the posets on $n$ points with at least $\tau \cdot 2^n$ downsets, {\em Order} {\bf 36} (2019), 119--157. Published Online May 29, 2018, \url{https://doi.org/10.1007/s11083-018-9459-2}.

\bibitem{aCampo_Erne_2019} F. a Campo and M. Ern\'{e}, Exponential functions of finite posets and the number of extensions with a fixed set of minimal points, {\em J.\ Combin.\ Math.\ Combin.\ Comput.} {\bf 110} (2019), 125--156.

\bibitem{Church_1940} R. Church, Numerical analysis of certain free distributive structures, {\em Duke Math.\ J.} {\bf 6} (1940), 732--734.

\bibitem{Church_1965} R. Church, Enumeration by rank of the elements of the free distributive lattice with 7 generators, {\em Notices Amer.\ Math.\ Soc.} {\bf 12} (1965), 724.

\bibitem{Davey_Priestley_2012} B. A. Davey and H. A. Priestley, {\em Introduction to Lattices and Order}, Cambridge University Press, 2\textsuperscript{nd} ed., 7\textsuperscript{th} printing 2012.

\bibitem{Dedekind_1897} R. Dedekind, \"{U}ber Zerlegungen von Zahlen durch ihre gr\"{o}ssten gemeinsamen Theiler, {\em Festschrift Hoch.\ Braunschweig u.\ ges.\ Werke II} (1897), 103--148.

\bibitem{Fidytek_etal_2001} R. Fidytek, A. W. Mostowski, R. Somla, and A. Szepietowski, Algorithms counting monotone Boolean functions, {\em Inform.\ Process.\ Lett.} {\bf 79} (2001), 203--209.

\bibitem{Lunnon_1971} F. Lunnon, The IU function: The size of a free distributive lattice, in D. J. A. Welsh ed., {\em Combinatorial Mathematics and its Applications}, Academic Press, 1971, pp.\ 173--181.

\bibitem{Markowsky_1989} G. Markowsky, Enumerating free distributive lattices, {\em Report University of Maine} {\bf 10} (1989).

\bibitem{Riviere_1968} N. M. Riviere, Recursive formulas on free distributive lattices, {\em J. Combin.\ Theory} {\bf 5} (1968), 229--234.

\bibitem{Shmulevich_etal_1995} I. Shmulevich, T. M. Sellke, M. Gabbouj, and E. J. Coyle, Stack filters and free distributive lattices: {\em Proceedings of 1995 IEEE Workshop on Nonlinear Signal Processing}, Halkidiki, Greece, 1995, pp.\ 927--930.

\bibitem{Sloane_DB} N. J. A. Sloane, {\em  On-line Encyclopedia of Integer Sequences}. The OEIS Foundation, \url{http://oeis.org/A000372}.

\bibitem{Ward_1946} M. Ward, Note on the order of the free distributive lattice, {\em Bull.\ Amer.\ Math.\ Soc.} {\bf 52} (1946), 423.

\bibitem{Wiedemann_1991} D. Wiedemann, A computation of the eighth Dedekind number, {\em Order} {\bf 8} (1991), 5--6.

\bibitem{Yusun_2011} T. J. Yusun, {\em Dedekind Numbers and Related Sequences}. Master Thesis, Simon Fraser University, London, 2011.

\end{thebibliography}
\end{document}